\documentstyle{amsppt}
\overfullrule=0pt
\font\sc=cmcsc10
\magnification=\magstep1
 \pagewidth{6.0 true in}
\pageheight{9.0 true in}

\hoffset 0in 
\voffset -0.10in

\parindent 0pt
\parskip 14pt
\hyphenation{e-qui-va-len-ti}
\topmatter
\title The Poincar\'e lemma and local embeddability
\endtitle
\author Judith Brinkschulte \\
C. Denson Hill \\
Mauro Nacinovich
\endauthor
\address Judith Brinkschulte: Universit\'e de Grenoble 1, Institut
Fourier, B.P. 74,
38402 St. Martin d'H\`eres, France\endaddress
\email brinksch\@ujf-grenoble.fr \endemail
\address C.Denson Hill - Department of Mathematics, SUNY at Stony Brook,
Stony Brook NY 11794, USA \endaddress
\email dhill\@math.sunysb.edu \endemail
\address Mauro Nacinovich - Dipartimento di Matematica "L.Tonelli" -
via F.Buonarroti, 2 - 56127 PISA, Italy \endaddress
\email nacinovi\@ dm.unipi.it\endemail
\abstract \par
\centerline{(Italian)}\par
Per variet\`a $CR$ astratte pseudoconvesse la validit\`a del Lemma di
Poincar\'e per forme di tipo $(0,1)$ implica
l'immergibilit\`a locale in $\Bbb C^N$; le due propriet\`a
sono equivalenti  per ipersuperfici
di dimensione reale $\geq 5$. Come corollario si ottiene un criterio
per la non validit\`a del Lemma di Poincar\'e per forme di tipo $(0,1)$
per una vasta classe di variet\`a $CR$ astratte di codimensione $CR$
maggiore di uno. \medskip
\centerline{(English)}\par
For pseudocomplex abstract $CR$ manifolds, the validity of the
Poincar\'e Lemma for $(0,1)$ forms implies local embeddability in
$\Bbb C^N$. The two properties are equivalent for hypersurfaces of
real dimension $\geq 5$. As a corollary we obtain a criterion for
the non validity of the Poicar\'e Lemma for $(0,1)$ forms for
a large class of abstract $CR$ manifolds of $CR$ codimension
larger than one.\endabstract
\endtopmatter

\document
\noindent
{\bf \S 1 Abstract $CR$ manifolds}
\smallskip

An abstract $CR$ manifold of type $(n,k)$ is a triple $(M, HM, J)$ where $M$
is a smooth real manifold of dimension $2n+k$, $HM$ is a subbundle of rank
$2n$ of the tangent bundle $TM$, and $J: HM \rightarrow HM$ is a smooth fiber
preserving bundle isomorphism with $J^2 =-I$. We also require that $J$ be
formally integrable; i.e. that we have
$$ [T^{0,1},T^{0,1}] \subset T^{0,1}$$
where
$$ T^{0,1} = \{ X+iJX \mid X \in \Gamma (M,HM)\} \subset \Gamma (M, \Bbb C
TM),$$
with $\Gamma$ denoting smooth sections.\par
We denote by $H^o M=\{\xi\in T^* M\mid <X,\xi >=0,\ \forall X\in H_{\pi
(\xi)}M\}$ the {\it characteristic bundle} of $M$. To each $\xi\in H^o_p M$,
we associate the Levi form at $\xi$:\newline
\centerline{${\Cal L}_p(\xi,X) = \xi ([J\widetilde{X},
\widetilde{X}])=
d\widetilde{\xi}(X,JX)$ for $X\in H_p M$}
which is
Hermitian for the complex structure of $H_p M$ defined by $J$.
Here
$\widetilde{\xi}$ is a section of $H^o M$ extending $\xi$ and
$\widetilde{X}$
a section of $HM$ extending $X$.\par
We denote by
$\bar\partial_M$ the tangential Cauchy-Riemann operator on $M$. A
smooth
function $f$ is called a $CR$ function on $M$ if
$\bar\partial_M
f=0$.\par
We say that $(M,HM,J)$ is locally $CR$ embeddable at
$p \in M$ if there exist
$n+k$ smooth complex valued $CR$ functions on a
neighborhood of $p$ whose
differentials are linearly independent.\par
We
say that the Poincar\'e lemma for $\bar\partial_M$ on $(0,1)$ forms
holds
at $p \in M$ if for every open neighborhood $\Omega$ of $p$ and every
smooth
$(0,1)$ form $f$ on $\Omega$ with $\bar\partial_M f=0$, there exists
an open
neighborhood $\omega \subset\Omega$ of $p$ and a smooth function
$u$ on
$\omega$ satisfying $\bar\partial_M u=f$ on $\omega$.
\medskip


\noindent
{\bf \S 2 Codimension one}\nopagebreak\smallskip\nopagebreak
Let $M$ be a smooth
abstract $CR$ manifold of type $(n,1)$ with $n \geq 2$,
and consider a
point $p \in M$.
\smallskip

\proclaim{Theorem 1} Assume $M$ is strictly
pseudoconvex at $p$. Then the
Poincar\'e lemma for $\bar\partial_M$ on
smooth forms of bidegree $(0,1)$ is
valid at $p$ if and only if $M$ is
locally $CR$ embeddable at $p$.
\endproclaim

\demo{Proof}
If $M$ is
locally $CR$ embeddable at $p$, then since $M$ is strictly
pseudoconvex and
$\dim_{\Bbb R}M \geq 5$, the Poincar\'e lemma for smooth
$(0,1)$ forms is
known to be valid by [AH].\par
For the proof in the other direction, we
shall employ the trick of Boutet de
Monvel [B] and use some a priori
estimates from [AFN].\par
Choose local coordinates $x_1,x_2, \ldots ,
x_{2n+1}$ for $M$, centered at
$p$, so that $p$ becomes the origin. By the
formal Cauchy-Kowalewski
procedure, we can find smooth complex valued
functions $ \varphi = (\varphi_1,
\varphi_2, \ldots ,\varphi_{n+1})$ in an
open neighborhood $U$ of $0$ with
each $\varphi_j (0) = 0$, $d\varphi_1
\wedge d\varphi_2
\wedge\ldots \wedge d\varphi_{n+1} \not= 0$ in $U$, and
such that
$\bar\partial_M \varphi_j$ vanishes to infinite order at $0$.
Then $\varphi : U
\longrightarrow \Bbb C^{n+1}$ gives a smooth local
embedding $\widetilde{M} =
\varphi (U)$ of $M$ into $\Bbb C^{n+1}$. On
$\widetilde{M}$ there is the $CR$
structure induced from $\Bbb C^{n+1}$; it
agrees to infinite order at $0$ with
the original $CR$ structure on $M$. In
particular $\widetilde{M}$ is a smooth
real hypersurface in $\Bbb C^{n+1}$
which is strictly pseudoconvex with
respect to the induced $CR$ structure.
This means that after possibly
shrinking $U$, there is a smooth complex
valued function $h$ in $U$ such that
$h(0)=0$, $\bar\partial_M h$ vanishes
to infinite order at $0$ and
$$ \Re h(x) \geq c \vert x\vert^2 \tag
\text{$*$}$$
with a positive constant $c$. Indeed it is well known that
after a suitable
local biholomorphic mapping, $\widetilde{M}$ can be
assumed to be strictly
convex with $T_o M=\{ \Re w_1 = 0 \}$ and $H_o M= \{
w_1 =0 \}$. It suffices
to take $h= \varphi^{*} w_1$.\par
Set
$w_j^{\lambda}= \bar\partial_M (\varphi_j e^{- \lambda h})$ for
$j=1,
\ldots ,n+1$ and $\lambda \in \Bbb N$. Then each $w_j^{\lambda}
\rightarrow 0$
in the topology of $\Cal C^{\infty}(U)$ as $\lambda
\rightarrow \infty$.
Indeed by ($*$) the function $\exp{\{-\lambda h\}}$,
and any derivative of it
with respect to $x$, is rapidly decreasing as
$\lambda \rightarrow \infty$,
while all other terms, and their derivatives
with respect to $x$, have only
polynomial growth in $\lambda$.\par By our
assumption that the Poincar\'e
lemma is valid at $0$, it follows from [AFN;
p. 384] that there are open
neighborhoods $W \subset V \subset\subset U$ of
$0$, a positive integer $l$
and a positive constant $A$ such that:\par For
every $\lambda \in \Bbb N$
there exists $u_j^{\lambda} \in {\Cal
C}^{\infty}(W)$ satisfying
$\bar\partial_M u_j^{\lambda} = w_j^{\lambda}$
in $W$ and
$$ \sup_{x\in W} \sup_{\vert \alpha \vert \leq 1} \vert
D^{\alpha}
u_j^{\alpha}\vert \leq A \sup_{x\in V} \sup_{\vert \beta\vert
\leq l} \Vert
D^{\beta}w_j^{\lambda}\Vert . \tag \text{$**$}$$
Here the
$\Vert \ \Vert$ on the right hand side indicates a norm on
$(0,1)$
forms.\par
Next we set $\psi_j^{\lambda}= \varphi_j e^{-\lambda h}-
u_j^{\lambda}$.
Then each  $\psi_j^{\lambda}$ is a smooth $CR$ function on
$W$. Since
$\varphi_j (0)=h(0)=0$, it follows that $\varphi_j$ and
$\varphi_j
\exp{\{-\lambda h\}}$ have the same first derivatives at $0$.
Finally we take
$\lambda$ sufficiently large. Then by ($**$) we have that
$d\psi_1^{\lambda}
\wedge d\psi_2^{\lambda} \wedge \ldots \wedge
d\psi_{n+1}^{\lambda}\not= 0$
on $W$. Thus $\psi = (\psi_1^{\lambda},
\ldots ,\psi_{n+1}^{\lambda}): W
\longrightarrow \Bbb C^{n+1}$ gives the
desired $CR$ embedding.
\enddemo
\smallskip
\noindent{\sc Remarks.}
In
Theorem 1 we have $\dim_{\Bbb R}M\geq 5$. When $\dim_{\Bbb R}M \geq 7$
the
local $CR$ embeddability is known to be possible (see [K], [A], [W]).
When
$\dim_{\Bbb R}M \geq 5$, and $M$ is $CR$ embedded, the Poincar\'e
lemma for
$(0,1)$ forms is also known to be valid (see [AH]). For abstract
strictly
pseudoconvex $M$ there are two open problems in dimension 5: the
local $CR$
embeddability, and the validity of the Poincar\'e lemma for
$(0,1)$ forms. By
Theorem 1, these two problems are equivalent. When
$\dim_{\Bbb R}M =3$, one
cannot always locally $CR$ embed by [Ni], and even
if $M$ is $CR$ embedded, the
Poincar\'e lemma for $(0,1)$ forms fails by
[AH].\par
We note that the two equivalent conditions in Theorem 1
are, in fact, equivalent to a third apparently weaker condition:
namely that the range of $\bar\partial_M$ on functions is 
"closed". By this we mean that given any open neighborhood $\Omega$ of $p$,
and any sequence $\{f_n\}_n$ of smooth functions on $\Omega$ such that
$\bar\partial_M f_n$ converges to $g$ in ${\Cal C}^{\infty}(\Omega )$,
there exists an open neighborhood $\omega$ of $p$ such that $g$ is
$\bar\partial_M$ exact on $\omega$. Indeed this condition is sufficient
to prove the local $CR$ embeddability at $p$, since one still has ($**$)
by [N].\par
Without strict pseudoconvexity, however, this "closed range property" is
no longer sufficient to obtain the local $CR$ embeddability. Indeed it
follows from the subelliptic estimates for functions on abstract $1$-
concave $CR$ manifolds proved in [HN1] that one has this closed range
property. However, there are examples of $1$-concave and even $2$-concave
$CR$ manifolds of type $(n,1)$ which cannot be locally $CR$ embedded (see
[JT], [R]).


\medskip
\noindent
{\bf \S 3 Higher codimension}
\smallskip
Let
$M$ be a smooth abstract $CR$ manifold of type $(n,k)$ with $n\geq 1$
and
$k\geq 1$, and consider a point $p \in M$.
\smallskip
\proclaim{Theorem 2}
Assume there is a $\xi\in H_p^o M$, $\xi \not=0$, with
${\Cal L}_p (\xi
,\cdot )$ positive definite. If the Poincar\'e lemma for
$\bar\partial_M$
on forms of bidegree $(0,1)$ is valid at $p$ then $M$ is
locally $CR$
embeddable at $p$.
\endproclaim
\demo{Proof} The proof is essentially the same as for Theorem 1, with $\Bbb
C^{n+1}$ replaced by $\Bbb C^{n+k}$, and with $j=1,2,\ldots ,n+k$. The first
crucial point is to have the estimate ($**$); fortunately we have it again
from [AFN; p. 384]. The second crucial point is to observe that the approximate
$CR$ embedding $\widetilde{M}$ in $\Bbb C^{n+k}$, which now has real
codimension $k$, is contained in a strictly pseudoconvex hypersurface. Indeed
for $\xi\in H_p^o M= H_0^o \widetilde{M}$ we have that the Levi form
$\widetilde{\Cal L}_0(\xi, \cdot )$ of $\widetilde{M}$ is positive definite,
and we may choose smooth real local defining functions $\rho_1, \rho_2, \ldots
,\rho_k$ for $\widetilde{M}$ near the origin in $\Bbb C^{n+k}$, with
$d\rho_1(0), d\rho_2(0), \ldots , d\rho_k(0)$ orthonormal with respect to the
standard Euclidean metric, such that locally
$$ \widetilde{M}= \{ \rho_j(z)=0 \mid j=1,\ldots ,k \}$$
and $d\rho_1(0)= \widetilde{J}^* \xi$, where $\widetilde{J}^* = -
\widetilde{J}$ is the adjoint of the complex structure tensor in $\Bbb
C^{n+k}$. Then replacing the $\rho_j$ by

\centerline{$\widetilde{\rho}_1 = \rho_1 \pm B \sum_{j=2}^k
\rho_j^2$}
 \centerline{$ \widetilde{\rho}_j = \rho_j$ for $2\leq j \leq k$}
\noindent with a large positive constant $B$, and choosing the appropriate
sign, we obtain a strictly pseudoconvex real hypersurface $\Sigma$, defined by
$\Sigma = \{ \widetilde{\rho}_1 (z)=0\}$, on which $\widetilde{M}$ lies.
As before this gives us the
existence of a smooth function $h$ such that $\bar\partial_M h$ vanishes to
infinite order at the origin, and satisfies ($*$), and the proof of Theorem 2
is complete. \enddemo

\proclaim{Corollary} Suppose\roster
\item "(a)" There exists $\xi \in H_p^o M$ with ${\Cal L}_p (\xi  ,\cdot )$
having one negative eigenvalue, and all other eigenvalues
greater than or equal to zero;
\item "{}"\quad \qquad and
\item "(b)" There
exists $\eta \in H_p^o M$ with ${\Cal L}_p (\eta  ,\cdot )$
positive definite.\endroster
\noindent Then the Poincar\'e lemma for $\bar\partial_M$ fails for $(0,1)$
forms at $p$. \endproclaim

\demo{Proof} Assume that the Poincar\'e lemma for $(0,1)$ forms is valid
at
$p$. Then by virtue of (b) and Theorem 2, $M$ is locally $CR$ embeddable
at
$p$. But then the Poincar\'e lemma for $(0,1)$ forms fails, using the
results
of [HN2]; this is a contradiction, completing the
proof.
\enddemo
\smallskip
\noindent{\sc Remarks.} Note that if $M$ is
assumed a priori to be $CR$
embedded at $p$, we have the result of the
Corollary dropping the hypothesis
(b). If in this situation the remaining
hypothesis (a) is strengthened to say
that all the other eigenvalues are
greater than zero, then the result already
follows from [AFN].\par
Note
that generically, when $k>1$, we expect
(a) to be a consequence of (b).
Indeed, if
${\Cal L}_p (\eta ,\cdot )$ has $n$ positive eigenvalues, then ${\Cal
L}_p
(-\eta ,\cdot )$ has $n$ negative eigenvalues. Hence if we move along
a
continuous path $\gamma (t)$, $0\leq t\leq 1$ from $\gamma (0)=\eta$
to
$\gamma (1)=-\eta$ on the sphere $S^{k-1}$, all eigenvalues of the Levi
form change sign. Since they are continuous functions of $t$,
there are points where the determinant of
the matrix of the Levi form vanishes. If at the first of these points,
say $t_0$,
the determinant of the Levi form has a simple zero, then
at a nearby
$\xi\in\gamma ([0,1])$, corresponding to a $t>t_0$
sufficiently close to $t_0$, (a) is satisfied.\par
Thus, generically we expect, by [HN2],
 that a locally embeddable $M$, of type
$(n,k)$ for some $k\geq 2$, and satisfying (b), does not admit
the Poincar\'e lemma for $(0,1)$ forms.
In the same situation, for an abstract $M$, the
Poincar\'e
lemma was known in general to fail for either $(0,1)$ forms or
for $(0,2)$
forms (see [N]).
Thus the Corollary above somewhat precises and improves, in the special case
of the tangential Cauchy-Riemann complex,
a
result of [N].\par

When $M$ is abstract and of type $(1,k)$ then
$\bar\partial_M u=f$ becomes a
scalar equation $Lu=f$, and the hypotheses (a) and (b)
in the Corollary are equivalent and both
can be restated by saying that $L,
\bar L , [L,\bar L ]$ are linearly
independent at $p$. Hence one obtains a
more geometric explanation of (a
special case of) H\"ormander's nonsolvability
result
[H\"o].

\enddocument